\documentclass{article}
\usepackage{amssymb,amsmath,amsthm,graphicx,CJK}

\def\qed{\hfill {\hbox{${\vcenter{\vbox{               
   \hrule height 0.4pt\hbox{\vrule width 0.4pt height 6pt
   \kern5pt\vrule width 0.4pt}\hrule height 0.4pt}}}$}}}

\def\utr{\, \underline{\ast}\, }
\def\otr{\, \overline{\ast}\, }

\theoremstyle{definition}
\newtheorem{example}{Example}
\newtheorem{definition}{Definition}

\date{}

\title{\Large \textbf{Virtual Links with Finite Medial Bikei}}

\author{
Julien Chien\footnote{Email: jchien17@cmc.edu} \and
Sam Nelson\footnote{Email: Sam.Nelson@cmc.edu. Partially supported by Simons Foundation collaboration grant 316709}}

\begin{document}
\maketitle

\begin{abstract}
We consider the question of which virtual knots have finite fundamental 
medial bikei. We describe and implement an algorithm for completing a 
presentation matrix of a medial bikei to an operation table, determining 
both the cardinality and isomorphism class of the fundamental 
medial bikei, each of which are link invariants. As an example, we compute 
the fundamental medial bikei for all of the prime virtual knots with up to 
four classical crossings as listed in the knot atlas.
\end{abstract}

\medskip

\parbox{5.5in}{\textsc{Keywords:} bikei, involutory biquandles, medial bikei, finite presentations}

\smallskip


\section{\large\textbf{Introduction}}\label{I}

In \cite{AN} an algebraic structure known as \textit{bikei} 
(\begin{CJK*}{UTF8}{min}双圭\end{CJK*}) was introduced, generalizing the notion 
of \textit{kei} (\begin{CJK*}{UTF8}{min}圭\end{CJK*}) or \textit{involutory 
quandles} from \cite{J,T}. Every unoriented 
classical or virtual knot or link has an associated \textit{fundamental bikei}
with the property that ambient isotopic knots and links have isomorphic
fundamental bikei. Moreover, various quotients of the fundamental bikei
obtained by imposing extra algebraic axioms have been defined, each also
invariant under ambient isotopy and hence defining invariants of knots and 
links. In particular, the notion of \textit{abelian} or \textit{medial} bikei
was defined in \cite{CN} for biquandles, which specialize to the case of bikei.

It is natural to ask whether some of these bikei or medial bikei are 
finite, and thus can be directly used to compare different knots.
However, determining which knots have finite bikeis and computing the structure
of these bikeis can take some work -- given a knot or link diagram, we 
can obtain
a \textit{presentation} of the fundamental bikei, but comparing isomorphism 
classes of objects defined by presentations can be nontrivial. 
In other work such as \cite{HS,NT} an algorithm with roots dating back to 
\cite{C} is used to compute the operation table of a finite algebraic 
structure given initially by a finite presentation.

In this paper, we describe and implement an algorithm for computing the 
fundamental medial bikei for virtual knots and links (including the 
classical case). In particular, we identify some
cases when these fundamental medial bikei of a virtual knot or
link is finite. The algorithm can be used to prove finiteness
and detect the isomorphism classes of fundamental medial bikei for some virtual knots and links. The cardinality of the fundamental medial bikei
is an integer-valued invariant of virtual knots and links when 
finite. For virtual knots and links with the same size medial
bikei, the isomorphism class is a generally stronger invariant than the
cardinality alone, since isomorphic bikei necessarily have equal cardinality. 
Moreover,
knowing the fundamental medial bikei structure can be helpful in determining 
which bikei to use for counting invariants and their enhancements, since
if two virtual links $L, L'$ have isomorphic fundamental medial bikei then 
for any finite medial bikei $X$ the counting invariants 
$\mathrm{Hom}(\mathcal{BK}(L),X)$ and $\mathrm{Hom}(\mathcal{BK}(L'),X)$ will
be equal, so we must then use non-medial bikei $X$ if we wish to detect 
any difference.

The paper is organized as follows. In Section \ref{V} 
we recall the basics of virtual knots and links. In section 
\ref{B} we give a brief review of bikei and medial bikei.
Section \ref{Com} contains a description of our algorithm and collects
computational results, 
and in Section \ref{Q} we collect some questions for future research.

\section{\large\textbf{Virtual Links}}\label{V}

We begin with a brief review of virtual knot theory. See \cite{K} for more.

A \textit{virtual link} is an equivalence class of \textit{virtual
link diagrams}, planar 4-valent graphs with vertices decorated as either 
\textit{classical crossings} or \textit{virtual crossings} as shown
\[\includegraphics{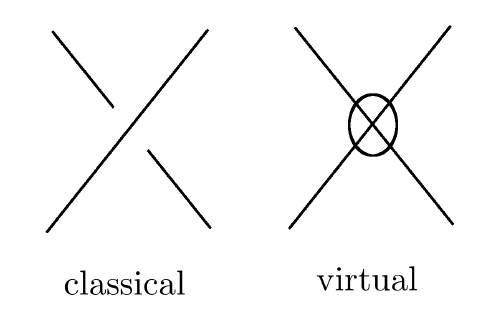}\]
under the equivalence relation generated by the seven 
\textit{virtual Reidemeister moves}:
\[\includegraphics{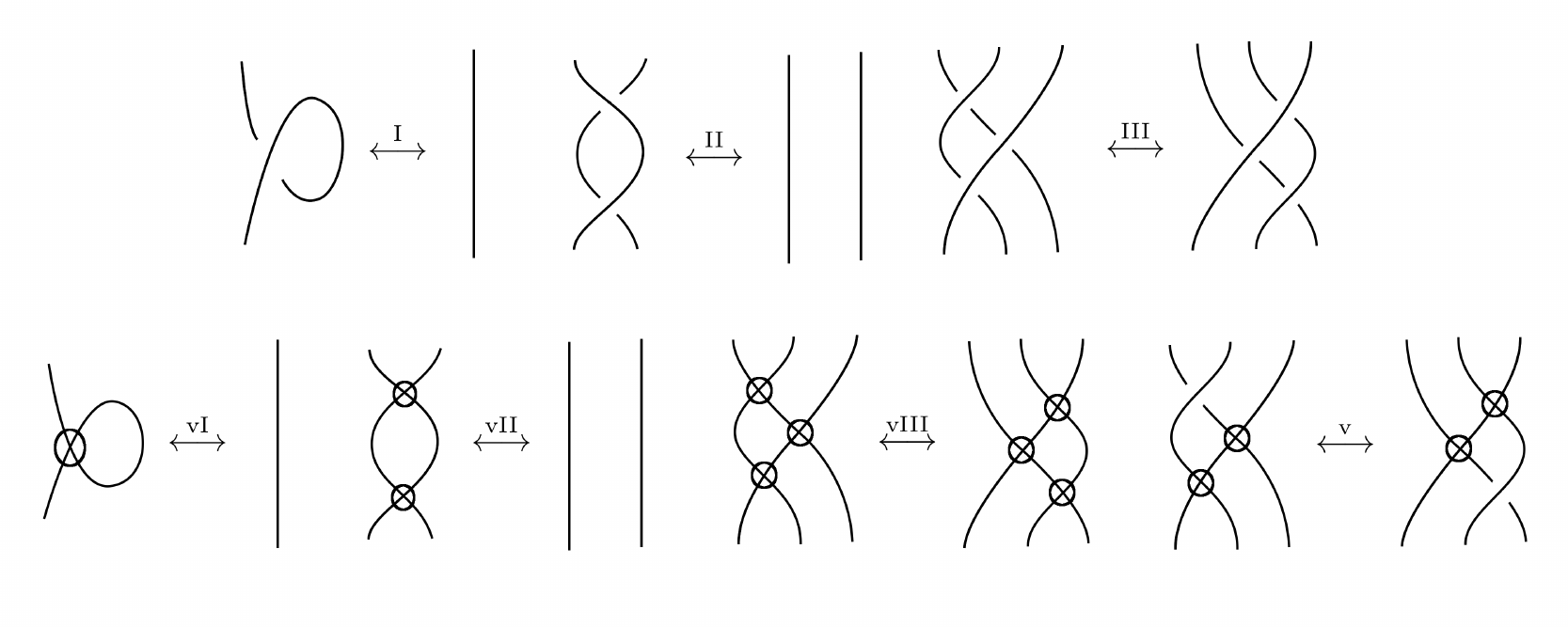}\]
In these moves, the knots or links in question are identical outside the
pictured portion of the diagrams.
A virtual link with a single component is a \textit{virtual knot}.

Virtual links can be interpreted as disjoint unions of simple closed curves 
in \textit{thickened surfaces}, i.e. orientable 3-manifolds-with-boundary of 
the form $\sigma\times [0,1]$ where $\sigma$ is a compact orientable surface, 
up to stabilization moves on $\sigma$. That is, we can think of drawing the 
virtual link diagram on a surface with possibly nonzero genus, with each virtual crossing 
representing a bridge or handle in the surface $\sigma$ and the classical
crossings drawn on $\sigma$; then, thickening $\sigma$, the classical crossings
represent points where the strands of the virtual link are close together inside
the thickened surface, while the virtual crossings result from compressing 
$\sigma$ onto genus-zero paper.
\[\includegraphics{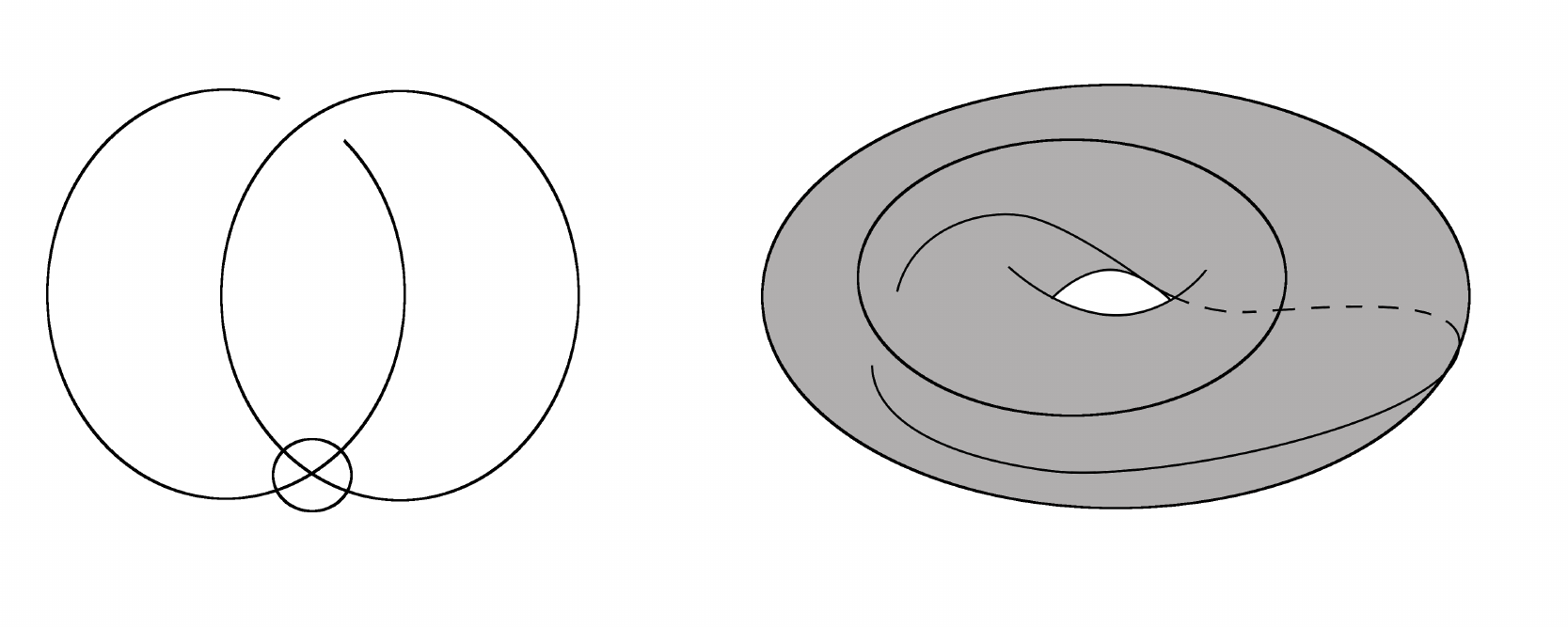}\]
A virtual knot is \textit{classical} if it is equivalent to a diagram with no
virtual crossings; these correspond to knots in ordinary genus zero
three-dimensional space. It is known (see \cite{KG} for instance) that
if two virtual crossing-free diagrams are equivalent by moves 
$\{\mathrm{I,\ II,\ III,\ vI,\ vII,\ vIII,\ v}\}$ then they
are also equivalent via moves $\{\mathrm{I,\ II,\ III}\}$ only.

Perhaps the main question in virtual knot theory is determining when two 
virtual knots or links are equivalent. This is generally done via 
\textit{virtual link invariants}, functions we can compute from virtual link
diagrams whose value does not change when the diagram is changed by  
Reidemeister moves. In the remainder of the paper we will describe our method
for computing two invariants, the size and isomorphism class of the 
fundamental medial bikei of a virtual link.

\section{\large\textbf{Bikei}}\label{B}

We begin this section with a definition (see \cite{AN,EN}).

\begin{definition}
A \textit{bikei} is a set $X$ with two binary operations
$\utr,\otr:X\times X\to X$ satisfying for all $x,y,z\in X$
\begin{itemize}
\item[(i)] $x\utr x=x\otr x$,
\item[(ii)]
\[\begin{array}{rclc}
x\utr(y\otr x)  & = & x\utr y & (ii.i)\\
x\otr(y\utr x)  & = & x\otr y & (ii.ii)\\
(x\otr y)\otr y & = & x & (ii.iii) \\
(x\utr y)\utr y & = & x & (ii.iv)\\
\end{array}\] and
\item[(iii)] (Exchange Laws)
\[\begin{array}{rclc}
(x\otr y)\otr (z\utr y) & = & (x\otr z)\otr (y\otr z) & (iii.i) \\
(x\otr y)\utr (z\otr y) & = & (x\utr z)\otr (y\utr z)  & (iii.ii)\\
(x\utr y)\utr (z\otr y) & = & (x\utr z)\utr (y\utr z) &  (iii.iii).
\end{array}\]
\end{itemize}
\end{definition}

\begin{example}
The integers $\mathbb{Z}$ and integers mod $n$ $\mathbb{Z}_n$
form bikei with operations $x\utr y=2y-x$ and $x\otr y=x$.
\end{example}

\begin{example}
A module $X$ over the ring $\mathbb{Z}/(t^2-1,s^2-1,(1-s)(s-t))$ is
a bikei under the operations $x\utr y=tx+(s-t)y$ and $x\otr y=sx$, known as
an \textit{Alexander bikei}. As a special case, $\mathbb{Z}_n$ is a bikei
under $x\utr y=tx+(s-t)y$ and $x\otr y=sx$ where we choose
$s,t\in \mathbb{Z}_n$ such that $s^2=t^2=1$ and $(1-s)(s-t)=0$.
\end{example}

Given a finite bikei $X$, we can represent the bikei structure with a block
matrix encoding the operation tables of $\utr,\otr$ in the following way:
let $X=\{x_1,\dots, x_n\}$. Then the \textit{bikei matrix} of $X$ is the
$n\times 2n$ matrix whose entry in row $j$ column $k$ is given by
$l\in\{1,\dots, n\}$ where
\[x_l=\left\{\begin{array}{ll}
x_j\utr x_k & 1\le k\le n \\
x_j\otr x_{k+n} & n+1\le k\le 2n \\
\end{array}\right.\]

\begin{example}
Let $X=\mathbb{Z}_4$ and set $t=1$ and $s=3$; then we have $s^2=t^2=1$ and
$(1-s)(s-t)=(1-3)(3-1)=2(2)=0$, so we have a bikei with operations
\[x\utr y=tx+(s-t)y=x+2y\quad\mathrm{and}\quad x\otr y=sx=3x.\]
Then using $X=\{1,2,3,4\}$ with the class of $0\in\mathbb{Z}_4$ represented
by $4$ so we can start our row and column numbering with 1, $X$ has
matrix
\[\left[\begin{array}{rrrr|rrrr}
3 & 1 & 3 & 1 & 3 & 3 & 3 & 3 \\
4 & 2 & 4 & 2 & 2 & 2 & 2 & 2 \\
1 & 3 & 1 & 3 & 1 & 1 & 1 & 1 \\
2 & 4 & 2 & 4 & 4 & 4 & 4 & 4
\end{array}\right].\]
\end{example}

\begin{example}
If $X$ is a bikei in which $x\otr y=x$ for all $X$, then $X$ is
called a \textit{kei} or  \textit{involutory quandle}.  For example,
any group $G$ forms a kei called a \textit{core kei} with 
$x\utr y= yx^{-1}y$ and $x\otr y=x$; the group 
$S_3=\{x_1=(1),x_2=(12),x_3=(13),x_4=(23),x_5=(123),x_6=(132)\}$
has bikei matrix
\[\left[\begin{array}{rrrrrr|rrrrrr}
1 & 1 & 1 & 1 & 6 & 5 & 1 & 1 & 1 & 1 & 1 & 1 \\
2 & 2 & 4 & 3 & 2 & 2 & 2 & 2 & 2 & 2 & 2 & 2 \\
3 & 4 & 3 & 2 & 3 & 3 & 3 & 3 & 3 & 3 & 3 & 3 \\
4 & 3 & 2 & 4 & 4 & 4 & 4 & 4 & 4 & 4 & 4 & 4 \\
5 & 6 & 6 & 6 & 5 & 1 & 5 & 5 & 5 & 5 & 5 & 5 \\ 
6 & 5 & 5 & 5 & 1 & 6 & 6 & 6 & 6 & 6 & 6 & 6
\end{array}\right]\]
\end{example}

\begin{definition}
A bikei $X$ is \textit{medial} or \textit{abelian} if it satisfies
\[
\begin{array}{rcll}
(x\utr y)\utr (z\utr w) & = &(x\utr z)\utr (y\utr w) & (m.i)\\
(x\utr y)\otr (z\utr w) & = &(x\otr z)\utr (y\otr w) & (m.ii)\\
(x\otr y)\otr (z\otr w) & = &(x\otr z)\otr (y\otr w) & (m.iii)\\
\end{array}
\]
for all $x,y,z,w\in X.$
\end{definition}

\begin{example}
Alexander bikei are always medial, since we have
\begin{eqnarray*}
(x\utr y)\utr (z\utr w) 
& = & t(tx+(s-t)y)+(s-y)(tz+(s-t)w)\\
& = & t^2x+t(s-t)(y+z)+(s-t)^2w\\
& = & t(tx+(s-t)z)+(s-y)(ty+(s-t)w)\\
& = & (x\utr z)\utr (y\utr w)
\end{eqnarray*}
so (m.i) is satisfied,
\begin{eqnarray*}
(x\utr y)\otr (z\utr w) 
& = & s(tx+(s-t)y)\\
& = & t(sx)+(s-t)(sy) \\
& = & (x\otr z)\utr (y\otr w)
\end{eqnarray*}
so (m.ii) is satisfied, and
\[(x\otr y)\otr (z\otr w)=s^2x=(x\otr z)\otr (y\otr w) \]
so (m.iii) is satisfied.
\end{example}

\begin{example}
The core bikei of a group $G$ need not be medial if $G$ is non-abelian,
since (m.i) requires
\[(x\utr y)\utr (z\utr w) =wz^{-1}wy^{-1}xy^{-1}wz^{-1}w
=wy^{-1}wz^{-1}xz^{-1}wy^{-1}w=
(x\utr z)\utr (y\utr w).\]
Indeed, for $G=S_3$ consider $x=(12)$, $y=(13)$, $z=(23)$ and $w=1$; then
\[(x\utr y)\utr (z\utr w) =wz^{-1}wy^{-1}xy^{-1}wz^{-1}w=(23)(13)(12)(13)(23)
=(23)
\]
while
\[(x\utr z)\utr (y\utr w)=wy^{-1}wz^{-1}xz^{-1}wy^{-1}w
=(13)(23)(12)(13)(23)=(12)\ne(23).
\]
\end{example}

Let $D$ be an unoriented virtual link diagram representing a virtual link
$L$. The \textit{fundamental bikei} of $L$ is the set of equivalence classes 
of \textit{bikei words} in a set $X$ of generators corresponding one to one 
with the \textit{semiarcs} of $D$, i.e., the portions of $D$ between crossing 
points, modulo the equivalence relation generated by the
bikei axioms and the \textit{crossing relations} of $D$. More precisely, 
let $X=\{x_1,\dots,x_n\}$ be a set of generators, one for each semiarc in
$D$. Then define a set $W(X)$ recursively by the rules
\begin{itemize}
\item[(i)] $x\in X$ implies $x\in W(X)$ and
\item[(ii)] $x,y\in W(X)$ implies $x\otr y, x\utr y\in W(X)$.
\end{itemize}
Note that since the operations $\otr,\utr$ are not associative, we need 
parentheses, e.g., $x_1\in W(X)$ and $x_2\utr x_3\in W(X)$ implies 
$x_1\otr (x_2\utr x_3)\in W(X)$ etc.
The \textit{free bikei} on $X$ is the set of equivalence classes of elements
of $W(X)$ under the equivalence relation generated by the bikei axioms; that is,
we have $x\otr x\sim x\utr x$, $(x\otr y)\otr y\sim x$, etc.

Then the fundamental bikei $\mathcal{BK}(L)$ is the set of equivalence classes
of bikei words under the stronger equivalence relation generated by both the 
bikei axioms and the \textit{crossing relations}, namely
\[\includegraphics{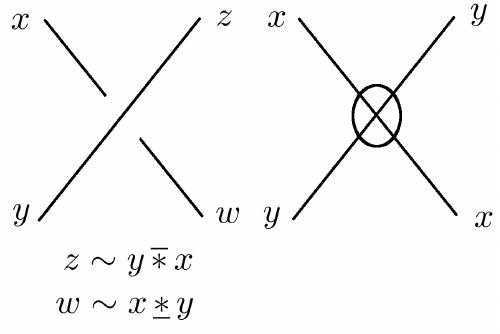}\]

We can specify such a bikei with a \textit{presentation} 
$\langle x_1,\dots, x_n \ |\ r_1,\dots r_m\rangle$ listing
generators $x_1,\dots, x_n$ and relations $r_1,\dots, r_m$. These relations are 
the equivalences generating the equivalence relation other than the bikei
axioms, which we don't list explicitly. Two such presentations describe 
isomorphic bikei if and only if they are related by \textit{Tietze moves},
which come in two types:
\begin{itemize}
\item[(i)] adding or removing a generator $x$ and relation of the form 
$x\sim W$ where $W$ is a word in the other generators not involving $x$, and
\item[(ii)] adding or removing a consequence of the other relations.
\end{itemize}

\begin{example}\label{ex:vt}
The pictured \textit{virtual trefoil} has fundamental bikei with four 
generators $x_1,x_2,x_3,x_4$ and relations 
$x_2=x_1\utr x_3,$ 
$x_4=x_3\otr x_1,$
$x_3=x_1\otr x_2$ and 
$x_4=x_2\utr x_1$
\[\includegraphics{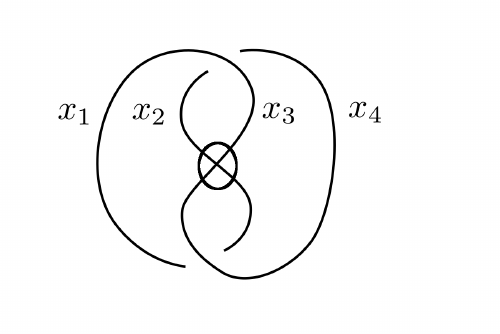}\]
\[\mathcal{B}(K)=\langle x_1,x_2,x_3,x_4\ |\ 
x_2=x_1\utr x_3,\ 
x_4=x_3\otr x_1,\ 
x_3=x_1\otr x_2,\ 
x_4=x_2\utr x_1 \rangle\]
\end{example}

Via Tietze moves, every bikei presentation can be put into 
\textit{short form}, in which every relation has the form $x_j=x_k\utr x_n$ 
or $x_j=x_k\otr x_n$ where $x_j,x_k$ and $x_n$ are generators as opposed to 
longer words. Given a presentation in short form, we can write the presentation
in matrix format with zeroes in positions without corresponding relations.

\begin{example}
Consider the bikei presented by 
$\langle x_1,x_2\ |\ (x_1\otr x_2)\utr x_1 = x_2\otr x_1\rangle$.
We can put this in short form by introducing new generators and corresponding 
relations. Via a Tietze I move, 
\[\langle x_1,x_2\ |\ (x_1\otr x_2)\utr x_1 = x_2\otr x_2\rangle\]
becomes
\[\langle x_1,x_2, x_3\ |\ x_3=x_1\otr x_2,\ x_3\utr x_1 = x_2\otr x_2\rangle\]
which then becomes
\[\langle x_1,x_2, x_3,x_4\ |\ x_3=x_1\otr x_2,\ x_3\utr x_1 = x_4, \
x_4=x_2\otr x_2\rangle.\]
Then in matrix form, we have presentation
\[\left[\begin{array}{rrrr|rrrr}
0 & 0 & 0 & 0 & 0 & 3 & 0 & 0 \\
0 & 0 & 0 & 0 & 0 & 4 & 0 & 0 \\
4 & 0 & 0 & 0 & 0 & 0 & 0 & 0 \\
0 & 0 & 0 & 0 & 0 & 0 & 0 & 0 \\
\end{array}\right].\]
\end{example}

\section{\large\textbf{Computation of the Fundamental Medial Bikei}}\label{Com}

If a bikei $X$ is finite, then its operation matrix is a short form
presentation matrix. Hence, we can play the following game: start with a 
medial bikei presentation, put it in short form and write its matrix, then 
use the medial bikei axioms to fill in entries by adding relations which are
consequences of the other relations, i.e. via Tietze II moves, in the hopes
of completing the table. If the table cannot be completed, we can select 
a zero entry in the table, assign to it a new generator and relation (a
Tietze I move), a new row and column in each matrix corresponding to the new generator,
and repeat the process of trying to fill in entries. 

During this process of filling in the matrices, we may encounter situations
where we need to assign a number to an entry in the matrix which is already
nonzero; if the two values for the position are equal, then there is nothing 
to do, but if they are distinct, then we have found that the two generators
assigned to that position are equivalent. We can then systematically replace all
instances of one generator with the other. Thus, we can loop over 
sets of generators, filling in entries using the medial bikei axioms, keeping a working list of which generators are equal, with the matrix growing or 
shrinking as new generators are consolidated or added. Since there is no
guarantee that the presented object is finite, this process may not terminate, 
but if it does, we have a proof that the object is finite in the form of its 
complete operation tables. To address the issue of non-terminating searches, we 
include a parameter defining the maximal size of matrix and terminate the 
computation when the matrices reach this size.

In terms of pseudocode:
\begin{verbatim}
//inputting the medial bikei presentation MB and maximal size L into the function:
Bikei(MB, L)
  M = shortFormMatrix(MB) // convert the medial bikei to matrix form
  I = 0
  Do while (M still has entries that equal 0) and (size(M) < L)
    Do while (ChangedM != M) 
      M = ChangedM
      sameGenerators = {} // keep track of equivalent generators
      //
      Do while (I <= size(M))
        ChangedM = firstBikeiAxiom(M,I,sameGenerators)
        J = 0
        Do while (J <= size(M))
          ChangedM = secondBikeiAxiom(ChangedM,I,J,sameGenerators)        
          Z = 0
          Do while (Z <= size(M))
            ChangedM = thirdBikeiAxiom(ChangedM,I,J,Z,sameGenerators)
            W = 0
            Do while (W <= size(M))
              ChangedM = medialBikeiAxiom(ChangedM,I,J,Z,W,sameGenerators)
            End-Do
          End-Do
        End-Do
      End-Do
      If (size(sameGenerators) > 0)
        ChangedM = reduce(ChangedM, sameGenerators) // merge equivalent entries
      End-If
    End-Do
    M = AssignNewGenerators(M) // assigns new generator to zero entries and expand M
  End-Do
  return M
\end{verbatim}

Our implementation in \texttt{Python} is available at \texttt{www.esotericka.org}.

\begin{example}
\[\includegraphics{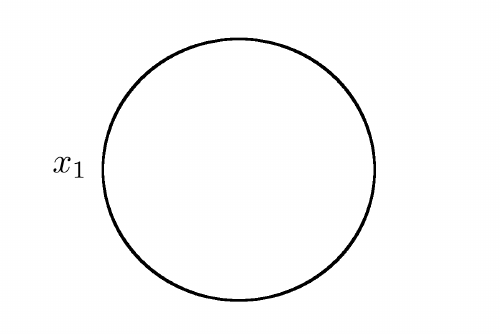}\]
The unknot at first may seem to have infinite fundamental medial bikei
since it has presentation $\langle x_1\ |\ \rangle$ with a single generator
and an empty list of relations. However, writing this in matrix form, we have
\[\left[\begin{array}{r|r}
0 & 0 \end{array}\right],\]
and applying our procedure, we add a new generator $x_2=x_1\utr x_1$:
\[\left[\begin{array}{rr|rr}
2 & 0 & 0 & 0\\
0 & 0 & 0 & 0 \end{array}\right].\]
Then axiom $(i)$ says $x\utr x=x\otr x$ and gives us
\[\left[\begin{array}{rr|rr}
2 & 0 & 2 & 0\\
0 & 0 & 0 & 0 \end{array}\right]\]
and axioms (ii.iii) and (ii.iv),  i.e., $(x\otr y)\otr y=x$ and
$(x\utr y)\utr y=x$ give us
\[\left[\begin{array}{rr|rr}
2 & 0 & 2 & 0\\
1 & 0 & 1 & 0 \end{array}\right].\]
Then axiom (ii.i), $x\utr (y\otr x)=x\utr y$, says
\[1\utr 2=1\utr(2\otr 1)=1\utr 1=2\]
and we have
\[\left[\begin{array}{rr|rr}
2 & 2 & 2 & 0\\
1 & 0 & 1 & 0 \end{array}\right];\]
similarly axiom (ii.ii), $x\otr(y\utr x)=x\otr y$
\[1\otr 2=1\otr(2\utr 1)=1\otr 1=2\]
and our matrix is
\[\left[\begin{array}{rr|rr}
2 & 2 & 2 & 2\\
1 & 0 & 1 & 0 \end{array}\right].\]
Lastly, another application of axioms (ii.iii) and (ii.iv), 
i.e., $(x\otr y)\otr y=x$ and $(x\utr y)\utr y=x$, yields
\[\left[\begin{array}{rr|rr}
2 & 2 & 2 & 2\\
1 & 1 & 1 & 1 \end{array}\right]\]
which says the fundamental medial bikei (in fact, plain fundamental bikei 
since we've not used the medial condition) of the unknot is the two element
bikei given by $\mathbb{Z}_2$ with $x\utr y=x\otr y=x+1$.
We note that this presentation can be read directly from the unknot diagram
below:
\[\includegraphics{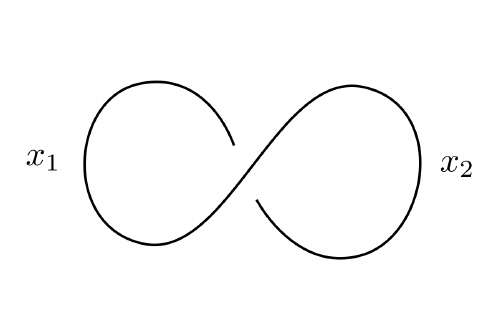}\]
\end{example}

\begin{example}
Consider the diagram below of the unknot:
\[\includegraphics{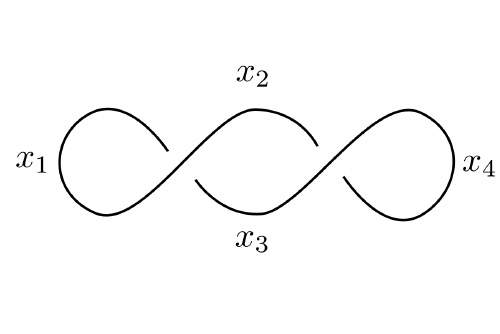}\]
We obtain the presentation matrix
\[\left[\begin{array}{rrrr|rrrr}
3 & 0 & 0 & 0 & 2 & 0 & 0 & 0 \\
0 & 0 & 4 & 0 & 0 & 0 & 0 & 0 \\
0 & 0 & 0 & 0 & 0 & 4 & 0 & 0 \\
0 & 0 & 0 & 0 & 0 & 0 & 0 & 0 \\
\end{array}\right].\]
The axiom (i) now says that $x_3=x_2$ since we have $x_1\utr x_1=x_3$ and 
$x_1\otr x_1=x_2$; then we must replace every instance of $3$ with $2$ and 
replace each instance of $j>3$ (as row number, column number, and entry)
with $j-1$. Note that in general 
this entails merging columns $2$ and $3$ and rows $2$ and $3$ in both 
matrices, which can trigger further merging. In this case we have
\[\left[\begin{array}{rrr|rrr}
2 & 0 & 0  & 2 & 0 &  0 \\
0 & 3 & 0  & 0 & 3 &  0 \\
0 & 0 & 0  & 0 & 0 &  0 \\
\end{array}\right].\]
Then axioms (ii.iii) and (ii.iv) give us
\[\left[\begin{array}{rrr|rrr}
2 & 0 & 0  & 2 & 0 &  0 \\
1 & 3 & 0  & 1 & 3 &  0 \\
0 & 2 & 0  & 0 & 2 &  0 \\
\end{array}\right].\]
Then as before, (ii.i) says
\[1\utr 2=1\utr(2\otr 1)=1\utr 1=2\]
so we have
\[\left[\begin{array}{rrr|rrr}
2 & 2 & 0  & 2 & 0 &  0 \\
1 & 3 & 0  & 1 & 3 &  0 \\
0 & 2 & 0  & 0 & 2 &  0 \\
\end{array}\right]\]
but then (ii.iii) says we must have $x_2=x_3$; reducing, we have
\[\left[\begin{array}{rr|rr}
2 & 2  & 2 & 2  \\
1 & 1  & 1 & 1  \\
\end{array}\right]\]
as before.
\end{example}

\begin{example}
The virtual trefoil knot in example \ref{ex:vt} has presentation matrix
\[\left[\begin{array}{rrrr|rrrr}
0 & 0 & 2 & 0 & 0 & 3 & 0 & 0 \\
4 & 0 & 0 & 0 & 0 & 0 & 0 & 0 \\
0 & 0 & 0 & 0 & 4 & 0 & 0 & 0 \\
0 & 0 & 0 & 0 & 0 & 0 & 0 & 0 \\
\end{array}\right].\]
We can then fill in entries using the medial bikei axioms. For instance,
bikei axiom (ii.iv) says $(x\utr y)\utr y=x$, so since
$x_2\utr x_1=x_4$ this says $x_4\utr x_1=x_2$ and our matrix becomes
\[\left[\begin{array}{rrrr|rrrr}
0 & 0 & 2 & 0 & 0 & 3 & 0 & 0 \\
4 & 0 & 0 & 0 & 0 & 0 & 0 & 0 \\
0 & 0 & 0 & 0 & 4 & 0 & 0 & 0 \\
2 & 0 & 0 & 0 & 0 & 0 & 0 & 0 \\
\end{array}\right].\]
Repeating this with the other bikei axioms, we obtain the matrix
\[\left[\begin{array}{rrrr|rrrr}
2 & 2 & 2 & 2 & 3 & 3 & 3 & 3 \\
4 & 1 & 1 & 1 & 1 & 1 & 1 & 1 \\
1 & 1 & 1 & 1 & 4 & 1 & 4 & 1 \\
2 & 2 & 2 & 2 & 3 & 3 & 3 & 3 \\
\end{array}\right]\]
together with the requirement that $x_1=x_4$ and $x_2=x_3$, which collapses
to
\[\left[\begin{array}{rr|rr}
2 & 2 & 2 & 2 \\
1 & 1 & 1 & 1 \\
\end{array}\right].\]
Hence, the fundamental medial bikei of the virtual trefoil is the
same as that of the unknot, and the invariant in this case does not detect the
nontriviality of the virtual trefoil.
\end{example}

\begin{example}\label{ex:6}
The virtual knot $4.71$ below
\[\includegraphics{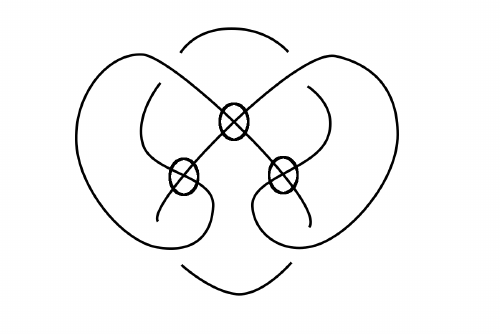}\]
has fundamental bikei presentation matrix
\[
\left[\begin{array}{rrrrrrrr|rrrrrrrr}
0& 0& 0& 0& 8& 0& 0& 0& 0& 0& 0& 2& 0& 0& 0& 0 \\
0& 0& 0& 0& 0& 0& 0& 0& 0& 0& 1& 0& 0& 3& 0& 0 \\
0& 4& 0& 0& 0& 0& 0& 0& 0& 0& 0& 0& 0& 0& 2& 0 \\
3& 0& 0& 0& 0& 0& 5& 0& 0& 0& 0& 0& 0& 0& 0& 0 \\
0& 0& 0& 0& 0& 0& 0& 4& 6& 0& 0& 0& 0& 0& 0& 0 \\
0& 7& 0& 0& 0& 0& 0& 0& 0& 0& 0& 0& 0& 0& 0& 5 \\
0& 0& 6& 0& 0& 0& 0& 0& 0& 0& 0& 8& 0& 0& 0& 0 \\
0& 0& 0& 0& 0& 1& 0& 0& 0& 0& 0& 0& 7& 0& 0& 0 \\ 
\end{array}\right]\]
which after application of our algorithm becomes
\[\left[\begin{array}{rrrrrr|rrrrrr}
5& 5& 5& 5& 5& 5& 5& 4& 2& 2& 5& 4 \\
6& 6& 6& 6& 6& 6& 3& 6& 1& 1& 3& 6 \\
4& 4& 4& 4& 4& 4& 2& 5& 4& 4& 2& 5 \\
3& 3& 3& 3& 3& 3& 6& 1& 3& 3& 6& 1 \\
1& 1& 1& 1& 1& 1& 1& 3& 6& 6& 1& 3 \\
2& 2& 2& 2& 2& 2& 4& 2& 5& 5& 4& 2 
\end{array}
\right].\]
This bikei is isomorphic to the Cartesian product of the unknot's bikei 
with the bikei obtained from the three-element Takasaki kei 
($\mathbb{Z}_3$ with $x\utr y=2x+2y$ and $x\otr y=x$) 
with
\[(x_1,x_2)\utr (y_1,y_2)=(x_1+1,x_2)\quad \mathrm{and}\quad
(x_1,x_2)\otr (y_1,y_2)=(x_1+1,2x_2+2y_2)
\]
after applying a vertical mirror image, i.e., switching $\otr$ and $\utr$.
Then since this is not the same bikei we obtained for the unknot, this
example shows the fundamental medial bikei distinguishing this virtual knot
from the unknot.
\end{example}

\begin{example}
We 
computed the fundamental medial bikei for all
virtual knots on the virtual knot table in \cite{KA}. Of the 116 virtual knots
on the list, most have the same medial fundamental bikei as the unknot computed above, but twenty-five do not; these break down into three isomorphism classes. 
According to our python computations.
\begin{itemize}
\item The virtual knots numbered $4.61$ through $4.77$ at \cite{KA}, up to 
vertical mirror image, all have fundamental medial bikei isomorphic to the 
six-element bikei listed in Example \ref{ex:6};
\item The virtual knots $3.6,\ 3.7,\ 4.98$ and $4.99$ all have fundamental 
medial bikei with 18 elements isomorphic to the fundamental medial bikei of the trefoil knot $3.6$, and
\item The virtual knots $4.105,4.106,4.107$ and $4.108$ all have fundamental 
medial bikei with 50 elements isomorphic to the fundamental medial bikei of the figure eight knot $4.108$.
\end{itemize}
\end{example}

\begin{example}
For examples with multicomponent virtual links, we do not have a convenient 
table analogous to the virtual knot table in \cite{KA}, but we note that 
our algorithm gives the eight-element medial bikei
\[
\left[\begin{array}{rrrrrrrr|rrrrrrrr}
2& 2& 1& 1& 2& 1& 2& 1& 2& 2& 5& 5& 2& 5& 2& 5 \\
1& 1& 2& 2& 1& 2& 1& 2& 1& 1& 7& 7& 1& 7& 1& 7 \\
6& 6& 4& 4& 6& 4& 6& 4& 3& 3& 4& 4& 3& 4& 3& 4 \\
8& 8& 3& 3& 8& 3& 8& 3& 4& 4& 3& 3& 4& 3& 4& 3 \\
7& 7& 5& 5& 7& 5& 7& 5& 7& 7& 1& 1& 7& 1& 7& 1 \\
3& 3& 8& 8& 3& 8& 3& 8& 6& 6& 8& 8& 6& 8& 6& 8 \\
5& 5& 7& 7& 5& 7& 5& 7& 5& 5& 2& 2& 5& 2& 5& 2 \\
4& 4& 6& 6& 4& 6& 4& 6& 8& 8& 6& 6& 8& 6& 8& 6 \\
\end{array}\right]
\]
for the virtual Hopf link, while the unlink of two components has
fundamental medial bikei given by the free medial bikei on two generators; 
this is infinite since it has the free Alexander kei on two generators
$\Lambda[x,y]$ where $\Lambda=\mathbb{Z}[t]/(t^2-1)$ as a quotient.
\end{example}

\section{\large\textbf{Questions}}\label{Q}

We conclude with a few questions for future research.

On the computational side, we have noticed that when choosing a zero entry 
to fill in with a new generator, the choice of which zero to fill in matters
greatly for the speed of completion of the procedure. We implemented a system
in which each zero receives a score based on the number of other entries 
which would be filled in as a result of filling in the zero in question; this 
seems to provide better results than simply selecting the zero based on 
position in the dictionary ordering, despite the associated computational 
overhead. However, in cases where the bikei is infinite, this may slow down
the procedure from reaching the exit size for the matrix.
What is the optimal strategy for zero selection?

On the mathematical side, the main question is what happens when we remove the 
medial condition or replace it with a different condition; which virtual links
have finite fundamental bikei?

\bibliography{jch-sn}{}
\bibliographystyle{abbrv}

\bigskip

\noindent
\textsc{Department of Mathematical Sciences \\
Claremont McKenna College \\
850 Columbia Ave. \\
Claremont, CA 91711}

\end{document}